\documentclass[12pt]{article}
\usepackage{latexsym,amssymb,amsmath,amsfonts,latexsym,amscd,amsthm}
\begin{document}

\newtheorem{ob}{Observation}[section]
    \newtheorem{lem}{Lemma}[section]
    \newtheorem{cor}{Corollary}[section]
    \newtheorem{theo}{Theorem}[section]
    \newtheorem{defi}{Definition}[section]
    \newtheorem{exam}{Example}[section]
    \newtheorem{rem}{Remark}[section]
	\newtheorem{prop}{Proposition}[section]
\numberwithin{equation}{section}
\title{
\Large
\bf
Convex Neighbourhoods and Complete Finsler Spaces
\footnote{
\newline Research
supported by  the Italian MURST  60\% and  GNSAGA}
}
 \author{
\small
O. M. Amici
and
  B. C. Casciaro\\
\small  Dipartimento di Matematica\\
\small  Universit{\`a} di Bari\\
\small Campus Universitario\\
\small Via Orabona 4, 70125 Bari, Italy
\\
\small amici@dm.uniba.it
\\
\small casciaro@dm.uniba.it
}
\maketitle
\begin{abstract}
In this paper, it is shown that a large set of
connections on a suitable sub--bundle of the tangent bundle
of a Finsler Manifold can be used to
study all the properties of convex neighbourhoods with respect to the Finsler Metric, which are needed to see that 
any Complete Finsler Space is Geodesically Connected.
\end{abstract}   
\newpage
\section{Introduction}
Let $M$ be a $C^\infty$--differentiable $n$--dimensional manifold endowed
with a Finslerian metric function $F:TM\to{\mathbb R}$, being $TM$ the tangent bundle
of $M$.

The following properties of the convex neighbourhoods of $M$
with respect to $F$ are
well known (see, \cite{BCS}, \cite{Ma},
see also
\cite{H}).

In order to quote them, let us denote by
\begin{eqnarray*}
B_\rho(0_x)=\{X\in T_x(M)/F(x,X)<\rho\}
\end{eqnarray*}
the
open indicatrix having $\rho$ as its radius and the zero vector
$0_x\in T_x(M)$ as its center, for any $\rho>0$ and $x\in M$.
\par
Moreover, we denote by $Exp_x$ the
exponential mapping of the Finsler metric and $Exp_x$ is defined on an open neighbourhood
of $0_x\in T_xM$ into $T_x(M)$, with $x\in M$. We recall that the mapping $Exp_x$ is defined by setting $Exp_x(X)=c_X(1)$, being $c_X$ the geodesic of $F$ defined by the initial conditions $c_X(0)=x$ and $\dot c_X(0)=X$, with $X$ belonging to a suitable open neighbourhood of $0_x$ in $T_xM$.
\par
With these notations, we have:
\begin{prop}
For each $x\in M$, there exist
two positive real numbers
$\varepsilon=\varepsilon(x)$ and $\eta=\eta(x)$ with $0<\varepsilon<\eta$
such that:

i). $Exp_x:B_\varepsilon(0_x)\to Exp_x(B_\varepsilon(0_x))=
B_\varepsilon(x)$ is a
diffeomorphism defined on $B_\varepsilon(0_x)$, whose degree of differentiability is
$C^\infty$ on $B_\varepsilon(0_x)-\{0_x\}$ and $C^1$ on $B_\varepsilon(0_x)$.
Moreover,
$B_\varepsilon(x)$ is an open neighbourhood of $x$.

ii). For each $y,z\in B_\varepsilon(x)$, there exists a unique geodesic
$c:[0,1]\to M$ lying entirely in
$B_\eta(y)=Exp_y(B_\eta(0_y))$ having length lesser than $\eta$ and such
that $c(0)=y$ and $c(1)=z$.

iii). For $y\in B_\varepsilon(x)$, the mapping $Exp_y:B_\eta(0_y)\to
B_\eta(y)$ is a diffeomorphism defined on $B_\eta(0_y)$, whose degree of
differentiability is $C^\infty$ on $B_\eta(0_y)-\{0_y\}$ and $C^1$ on
$B_\eta(0_y)$ and $B_\eta(y)$ is an open neighbourhood of
$y$ containing $B_\varepsilon(x)$.
\end{prop}

As in the Riemannian case (see, e.g., \cite{Kl}), we consider
the canonical identification
$T_XT_xM=T_xM$, for any $X\in T_xM$ and $x\in M$ and we fix an element $x$ of
$M$. Then, $T'_xM=T_xM-\{0_x\}$ is endowed by the Riemannian Metric
$g_p(X,Y)=g_{(x,p)}(X,Y)$, for each $p\in T'_xM$ and $X,Y\in T_pT_xM=T_xM$,
being $g$ the metric tensor which is induced by $F$.
Moreover, let $\widetilde TM\subseteq TM$ be the open neighbourhood of the zero
section where the exponential map of the Finsler Metric is defined.
Then, we have:
\begin{prop}
Let $x\in M$ and
$X\in \widetilde U_x=\widetilde TM\cap T_xM$. Moreover, let
$\widetilde b:[0,s_1]\to \widetilde U_x$ be
any differentiable curve such that $\widetilde b(0)=0_x$ and
$\widetilde b(s_1)=X$.
Finally let us put $b(s)=Exp_x\widetilde b(s)$, for any
$s\in[0,s_1]$ and $c(t)=Exp_x(tX)$, for any $t\in[0,1]$. Then:

i). $L(c)\leq L(b)$, where $L$ denotes the Finslerian length of
curves.

ii). If $\widetilde b(s)=\tau(s)X$, for any $s\in[0,s_1]$, being
$\tau:[0,s_1]\to[0,1]$ a strictly increasing differentiable mapping, then $L(c)=L(b)$.

iii). If $L(c)=L(b)$ and the total differential
$(DExp_x)_{t\widetilde b(s)}$ has maximal rank, for any $s\in[0,s_1]$ and $t\in[0,1]$, then
there exists a differentiable map $\tau:[0,s_1]\to[0,1]$ such that
$\widetilde b(s)=\tau(s)X$, for each $s\in[0,s_1]$ and $\tau$ is strictly increasing.
\end{prop}

\begin{prop}
Let $\widetilde \varepsilon=\varepsilon/3$, and $\varepsilon$ as in
Proposition 1.1. Then, every $y,z\in B_{\widetilde\varepsilon}(x)$
can be joined by a geodesic $c$,
lying entirely in
$B_{\widetilde\varepsilon}(x)$ and
having Finslerian length equal to the
Finslerian distance between $y$ and $z$. Moreover, any further geodesic joining $y$ and $z$ (if it exists) has points outside of $\overline B_\eta(y)$. 
\end{prop}

Finally, we recall the following definition: we say that $F$ is a {\it
Complete Finsler Metric}, if and only if the
distance which is induced by $F$ on $M$ is complete. Then it results:
\begin{prop}
Let $F$ be a Finslerian metric on $M$. Then the
following properties are equivalent:

i). $F$ is a complete Finslerian metric.

ii). There exists a point $p\in M$ such that any geodesic
starting from $p$ can be extended to the whole $\mathbb R$.

iii). Assertion ii) holds, for any $p\in M$.

Moreover, it results:

iv). If Assertion i) holds, any two points of $M$ can be joined by
a geodesic of $F$.
\end{prop}

In 1967, B. T. Hassan (see \cite{H})  proved Proposition 1.1
by using the
{\it Cartan
Connection} relative to $F$. Proofs of more or less complete
versions of Proposition 1.1
are given by several authors by means of connections different from
the Cartan one. In 1993, D. Bao and S. S. Chern proved the same
proposition by using a new connection, called the {\it Chern Connection}
(see \cite{BC}), which coincides with the {\it Rund
Connection} (see \cite{A}). 
\par
Finally, in \cite{BCS} the previous
proposition was proved by means of the {\it Cartan Form}.
The proofs of all the
previous propositions can be found in this book.
\par
Generally, all these proofs are different from the corresponding proofs of the
Riemannian case, more difficult than them and this fact explains the number of the
proofs of Proposition 1.1.
\par
Some questions arise for the previous observations. For example: Why
there exist so much proofs for Proposition 1.1? Why they are so different
from the proof needed in the Riemannian case? Does the used
connection play any role? 
\par
In \cite{AC} and \cite{AC1}, by using the bundle $T'M=TM-\sigma(M)$, where
$\sigma:M\to TM$ is the zero section, we determined all the connections by means of which
the Euler--Lagrange assumes the simplest form. These connections were called {\it
Finslerian Connections} and they are a large class of connections.
\par
In this paper we prove:
\begin{prop}
Any Finslerian Connection can be used to
prove all the previous propositions.
\end{prop}
We also show that all the proofs of the previous propositions can follow
in the closest way the corresponding proofs
used in the Riemannian case and this fact can be useful for the proof of
further results. We choose to follow \cite{Kl}.
\par
In our opinion, the problems previously listed raise because classical methods used in this context mix up the Finsler Geometry with the geometry of the tangent bundle of the manifold under consideration. 
\par
By using the previous considerations and the methods introduced here, in \cite{AC2} we determine all the properties corresponding to the ones considered here for a large class of Lagrangian Functions and in \cite{AC3} we show that the results of \cite{AC2} give a new sufficient condition for the geodesic connectedness of the generalized Bolza problem.

\section{Definitions and Proofs}

Let $M$ be a $C^\infty$--differentiable $n$--dimensional manifold.

The notations, which are more frequently used in the following, are:
\par

$\alpha$). $TM=\cup_{x\in M}T_xM$ is the tangent bundle of $M$ and
$\pi:TM\to M$ its natural projection. Moreover, for each $x\in M$, either
$X_x$ or
equivalently $(x,X)$ will denote the same element of $T_xM$ according to
this element is considered as a tangent vector at $x\in M$ or as a point of
$TM$;
hence $X_x=(x,X)$.
\par
$\beta$). ${\cal I}(M)=\bigoplus_{(r,s)\in{\mathbb N}^2}{\cal I}^r_s(M)$
is the algebra of
tensor fields on $M$,
where ${\cal I}^1_0(M)={\cal X}(M)$ and ${\cal I}^0_0(M)={\cal F}(M)$ are
the Lie algebra of vector fields and the ring of $C^\infty$--differentiable
real valued functions, respectively.
\par
$\gamma$). Let $\sigma:M\to TM$ be the zero section and we set
$T'M=TM-\sigma(M)$. $T'M$ is the open sub--bundle of $TM$ of non--zero tangent
vectors of $M$ and we shall denote by
$\pi':T'M\to M$ its canonical projection.
\par

$\delta$). ${\cal I}_{\pi'}=\bigoplus{\cal I}^r_{s\pi'}$is the ${\cal F}(T'M)$--module
(${\cal F}(M)$--module) of differentiable tensor fields along $\pi'$, with
${\cal I}^1_{0\pi'}={\cal X}_{\pi'}$ and ${\cal I}^0_{0\pi'}={\cal F}(T'M)$.
\par

$\theta$). If $(U,\varphi)$ is a local chart of $M$, we denote by
$(TU,T\varphi)$ the local chart canonically induced on $TM$, with
$TU=\pi^{-1}U$. We set
$\varphi=(x^1,\dots,x^n)$, $T\varphi=(x^1,\dots,x^n,X^1,\dots,X^n)$,
$e_i=\frac{\partial}{\partial x^i}$,
$\varepsilon_{\hat i}=\frac{\partial}{\partial X^{\hat i}},\ e^j=dx^j,\ \varepsilon^{\hat
j}=dX^{\hat j}$, for any $i,j,\hat i,\hat j\in\{1,\dots,n\}$.
\par

We recall that:

A {\it non--linear connection} $\nabla$ ({\it with three
indices}) can be regarded as an ${\mathbb R}$--bilinear mapping
$\nabla:{\cal X}_{\pi'}\times {\cal I}(M)\to {\cal I}_{\pi'}$ such that:
\par

i). $\nabla_{fX}(kY)=fX(k)Y+fk\nabla_XY,\ \forall X\in{\cal X}_{\pi'},\quad
\forall Y\in{\cal X}(M),\ \forall f,k\in{\cal F}(M)\ ;$

ii). $\nabla$ commutes with all the contractions;

iii). for each $X\in{\cal X}_{\pi'}$, with $X$ positively
homogeneous of degree $\rho$, $\nabla_XY$ is positively homogeneous of degree
$\rho$, for any $Y\in{\cal X}(M)$.

Moreover, $\nabla$ defines a horizontal lift $h:{\cal X}_{\pi'}\to{\cal X}(T'M)$
which can be extended in a trivial way to the whole tensor algebra ${\cal I}_{\pi'}$. This
extension is called {\it Matsumoto lift} (see, e.g., \cite{H}).
We also recall that the identity $id:T'M\to T'M$ can be considered as a
vector field along $\pi'$ ($id$ is the so--called {\it fundamental vector
field along $\pi'$}). The vector field $id^h$ will be said
to be {\it associated to $\nabla$} and one of its properties is:
\begin{prop}
$id^h$ can be extended to a spray defined on the
whole $M$, which will be denoted again with $id^h$. This extension is
$C^1$--differentiable on $TM$ and $C^\infty$--differentiable on $T'M$.
Moreover, a
curve $\gamma:[0,1]\to M$ is a non constant geodesic of $id^h$ if and only if
$\gamma$ is a path of $\nabla$.
\end{prop}

The spray $id^h$ defines an exponential map $Exp:\widetilde TM\to M$,
being $\widetilde TM$ an open neighbourhood of $\sigma(M)$. We have
$Exp(x,X)=c_{(x,X)}(1)$, being $c_{(x,X)}:[0,\varepsilon]\to M$, with $\varepsilon>1$, the
geodesic of $id^h$ having $(x,X)$ as initial condition, for all
$(x,X)\in\widetilde TM$.
The exponential map verifies the following properties:

(1). $Exp$ is $C^1$--differentiable on $\widetilde TM$ and
$C^\infty$--differentiable on $\widetilde TM-\sigma(M)$ (see \cite{L}, pg. 72).
\par

(2). $Exp$ has maximal rank in the zero vector $0_x\in T_xM$, for
any $x\in M$ (see \cite{GKM}, 2.8, Satz. (a), pg. 61).
\par

(3). The map $(\pi,Exp):\widetilde TM\to M\times M$ has maximal rank in
$0_x\in T_xM$, for any $x\in M$ (see \cite{GKM}, 2.8, Satz. (c), pg. 61).
\par

(4). The total differential of $Exp$, when it is restricted to $\widetilde TM\cap
T_xM$ and it is calculated in $0_x$, coincides with the identity map, for any
$x\in M$ (see \cite{L}, Th. 8, pg. 72).
\par

Moreover, as in the Riemannian case (see, e.g., \cite{Kl}) from
the previous properties, it follows:
\begin{prop}
There exists an open
neighbourhood $\widetilde W$ of $\sigma(M)$, such that
$Exp(\widetilde W)=W$ is open and:
\par
i). For each $y,z\in W$ there exists a unique geodesic
$c:[0,1]\to M$ such that $c(0)=y,\ c(1)=z$, $\ \dot c(0)\in \widetilde W$ and
$c(t)\in W$, for any $t\in[0,1]$.
\par

ii). For each $y\in W$, the map
$Exp_y:\widetilde W\cap T_yM\to Exp(\widetilde W\cap T_yM)=W(y)$
is a diffeomorphism of class $C^1$, of class $C^\infty$ on
$(\widetilde W\cap T_yM)-\{0_y\}$ and if $y\in W(x)$,
with $x\in M$, then $W(x)\subseteq W(y)$.
\par

iii). The mapping $(\pi,Exp)_{\vert\widetilde W}:\widetilde W\to W\times W$
is a
diffeomorphism of class $C^1$ and of class $C^\infty$ on
$\widetilde W-\sigma(M)$.
\end{prop}

Taking into account the previous proposition, the set $W(0_x)=\widetilde W\cap
T_xM$ is an open neighbourhood of $0_x$ in $T_xM$ and
$W(x)=Exp_x(W(0_x))$ is an open neighbourhood of $x$, for any $x\in M$.

Now, let $F:TM\to{\mathbb R}$
be a Finsler metric and $g$ its metric tensor. Then $g$ is
a family of Riemannian metrics on $M$ depending on $(x,X)\in T'M$ and such
that $g_{(x,X)}(X_x,X_x)=F^2(x,X)$, for any $x\in M$ and $(x,X)=X_x\in TM$.
Locally, we set $g=g_{ij}e^i\otimes e^j$, with
\begin{eqnarray*}
g_{ij}=\frac{1}{2}\frac{\partial^2 F^2}{\partial X^i\partial X^j}\ .
\end{eqnarray*}
Let $V=V(TM)$ be the vertical sub bundle of $TTM$; i.e., $V$ is the subset
of $TTM$ containing all the vectors tangent to the fibres of the tangent
bundle of $M$. We fix a further sub bundle $H$ of $TTM$ such that
$T_{(x,X)}TM=H_{(x,X)}\oplus V_{(x,X)}$, for any $(x,X)\in TM$. Two
canonical projections $P:TTM\to H$ and $Q:TTM\to V$ are associated to the
previous splitting into direct sum of $TTM$. These projections are
determined
by two tensor fields of type (1,1) denoted, by an abuse of notation, again by
$P$ and $Q$. Locally, $P$ and $Q$ are defined by:
\begin{eqnarray*}
P=\delta^i_je_i\otimes e^j-P^{\hat i}_j\varepsilon_{\hat i}\otimes e^j\
\hbox{and}\quad
Q=\delta^{\hat i}_{\hat j}\varepsilon_{\hat i}\otimes \varepsilon^{\hat j}+
P^{\hat i}_j\varepsilon_{\hat i}\otimes e^j\
.
\end{eqnarray*}
We recall that the distribution $H$ is said to be a
{\it non linear--connection with two
indices} if the functions $P^{\hat i}_j$ are positively homogeneous of degree
one (see, e.g.,\cite{Ma}).
Moreover, any connection on $M$ determines a couple of tensor fields
of the previous kind.
\par
Now, we turn to the general case and we set
\begin{eqnarray}
G=g^v=g_{ij}e^i\otimes e^j
\ ,
\end{eqnarray}
being
$v$ the extension to the tensor fields along $\pi'$ of the usual
vertical lift of tensors (see \cite{YI}).
\par
Then, from \cite{AC1}, it follows:
\begin{prop}
There exists a connection
$\nabla$ on $T'M$ such that
\par
i). $C^2_1(P\otimes T)=0$, where $T$ is the torsion tensor field
of $\nabla$ and $C^2_1$ is the contraction of the first contra--variant index
with the second covariant index of $P\otimes T$.
\par
ii). Let us denote by $\pi_*:TTM\to TM$ the total differential of
the canonical projection $\pi:TM\to M$ (see \cite{KN}).
Then, for each vector field
$X\in{\cal X}(T'(M))$, having $\pi_*(X_{(x,Z)})=Z_x$, for any $Z_x=(x,Z)\in
T'(M)$ it results:
\begin{eqnarray*}
(\nabla_XG)(X,Y)=(\nabla_YG)(X,X)=0\ ,\quad\forall Y\in{\cal X}(T'M)
\ .
\end{eqnarray*}

iii). The mapping $\nabla':{\cal X}_{\pi'}\times{\cal I}(M)\to{\cal I}_{\pi'}$
defined by
\begin{eqnarray*}
\nabla'_XY=\pi_*(\nabla_{X^c}Y^c)\ ,\quad\forall X,Y\in{\cal X}(M)
\end{eqnarray*}
being $c$ the complete lift (see \cite{YI}), is a non--linear connection with
three indices.
\end{prop}
A connection $\nabla$ verifying i), ii) and iii) of the previous
proposition is called {\it Finslerian Connection}. Moreover, the non--linear
connection $\nabla'$, defined by iii) of the same proposition, is called
{\it Berwald Connection deduced from $\nabla$}.

From \cite{AC1}, we also get:
\begin{prop}
Let us denote $\nabla$ by a Finslerian Connection
and
let $\nabla'$ be the Berwald Connection deduced from $\nabla$. Let
$Y=(\gamma,\dot\gamma):[a,b]\to T'M$ be a curve.
Then, $\gamma$ is a path of $\nabla'$ if and only if there exists a vector field
$Z\in{\cal X}(T'M)$ such that:

\begin{eqnarray*}
(\nabla_{\dot Y}{\dot Y})(t)=
Q_{(\gamma(t),\dot\gamma(t))}(Z_{(\gamma(t),\dot\gamma(t))})\ ,
\quad\forall
t\in[a,b]
\ .
\end{eqnarray*}
\end{prop}

In \cite{AC1}, it was also proved that the set of all the Finslerian Connections
can be
obtained in the following way.

Let $\widetilde\nabla$ be any connection on $T'M$ and $(U,\varphi)$
be a chart of
$M$. With respect to the chart $(T'U,T'\varphi)$, induced by $(U,\varphi)$
on $T'M$, we set:
\begin{eqnarray*}
& &
\nabla_{e_j}e_k
=
\widetilde\Gamma^{1i}_{jk}e_i+
\widetilde\Gamma^{5\hat i}_{jk}\varepsilon_{\hat i}
\ ,\
\nabla_{\varepsilon_{\hat j}}e_k
=\widetilde\Gamma^{2i}_{\hat jk}e_i+\widetilde\Gamma^{6\hat i}_{\hat jk}\varepsilon_{\hat i}
\nonumber
\\
& &
\nabla_{e_j}\varepsilon_{\hat k}
=\widetilde\Gamma^{3i}_{j\hat k}e_i+
\widetilde\Gamma^{7\hat i}_{j\hat k}\varepsilon_{\hat i}\ ,\
\nabla_{\varepsilon_{\hat j}}\varepsilon_{\hat k}
=\widetilde\Gamma^{4i}_{\hat j\hat k}e_i+
\widetilde\Gamma^{8\hat i}_{\hat j\hat k}\varepsilon_{\hat i}
\ .
\end{eqnarray*}
The $8n^3$ functions
$
(\widetilde\Gamma^{1i}_{jk},
\widetilde\Gamma^{2i}_{\hat jk},
\widetilde\Gamma^{3i}_{j\hat k},
\widetilde\Gamma^{4i}_{\hat j\hat k},
\widetilde\Gamma^{5\hat i}_{jk},
\widetilde\Gamma^{6\hat i}_{\hat jk},
\widetilde\Gamma^{7\hat i}_{j\hat k},
\widetilde\Gamma^{8\hat i}_{\hat j\hat k}
)$
are the local components of $\widetilde\nabla$.
By using the local components of $G$ and $P$ and the previous $8n^3$
functions, we obtain the following new $8n^3$ functions defined on $T'U$:
\begin{eqnarray*}
& &
2\Gamma^{1i}_{jk}=g^{im}(\partial_jg_{mk}+\partial_kg_{mj}-\partial_mg_{jk}+P^{\hat
h}_m\partial_{\hat h}g_{jk}-P^{\hat h}_j\partial_{\hat h}g_{mk}-P^{\hat h}_k\partial_{\hat
h}g_{mj})\ ,
\nonumber
\\
& &
\Gamma^{2i}_{\hat jk}=0\ ,\quad
\Gamma^{3i}_{j\hat k}=0\ ,\quad
\Gamma^{4i}_{\hat j\hat k}=0\ ,
\nonumber
\\
& &
\Gamma^{5\hat i}_{jk}=P^{\hat i}_t\widetilde\Gamma^{1t}_{jk}-P^{\hat i}_t\Gamma^{1t}_{jk}+
\widetilde\Gamma^{5\hat i}_{jk},
\nonumber
\\
& &
\Gamma^{6\hat i}_{\hat jk}=P^{\hat i}_t\widetilde\Gamma^{2t}_{\hat jk}+
\widetilde\Gamma^{6\hat i}_{\hat jk}\ ,\ \
\Gamma^{7\hat i}_{j\hat k}=
P^{\hat i}_t\widetilde\Gamma^{3t}_{j\hat k}+
\widetilde\Gamma^{7\hat i}_{j\hat k}\ ,\ \
\Gamma^{8\hat i}_{\hat j\hat k}=P^{\hat i}_t\widetilde\Gamma^{4t}_{\hat j\hat k}+
\widetilde\Gamma^{8\hat i}_{\hat j\hat k}
\ ,
\end{eqnarray*}
where we used the notations $\partial_if=\frac{\partial f}{\partial x^i}$
and $\partial_{\hat
i}f=\frac{\partial f}{\partial X^{\hat i}}$,
for any $f\in{\cal F}(T'M)$ and any $i,\hat
i\in\{1,\dots,n\}$.

The above functions are the local components of a connection $\nabla$ of
$T'M$. Since $\nabla$ verifies i), ii) and iii) of Proposition 2.3, it is a
Finslerian Connection. The connection $\nabla$ can be used in order to
determine all the Finslerian Connections. In fact, let $\nabla^1$ be a further
connection and let us denote by $N\in{\cal I}^1_2(T'M)$ the tensor field defined by
setting:
\begin{eqnarray*}
N(X,Y)=\nabla_XY-\nabla^1_XY\ ,\quad\forall X,Y\in{\cal X}(T'M)
\ .
\end{eqnarray*}
There exists a unique tensor field, $\widetilde N$, of type $(1,2)$ along $\pi'$,
such that:
\begin{eqnarray*}
\widetilde N(X,Y)=\pi_*N(X^c,Y^c)\ ,\quad\forall X,Y\in{\cal X}(M)
\end{eqnarray*}
$\widetilde N$ is called the {\it tensor field associated to $N$}. Now, we
denote by $d'$ the so called "derivation along the fibres",
which is defined by setting $d'f=(\partial_{\hat i}f)e^{\hat i}$, for any
$f\in{\cal F}(T'M)$ and is extended to the whole tensor algebra ${\cal I}_{\pi'}$ in
the well--known way.
Then, we have:
\begin{prop}
Under the previous assumptions, the connection
$\nabla^1$ is a Finslerian Connection, if and only if the following
conditions hold for $N$:

i). $\widetilde N_{(x,Z)}(Z_x,Z_x,)=0\ ,\quad\forall (x,Z)=Z_x\in T'M$.

ii). $C^1_1(Z_x\otimes C^1_1((d'F\otimes\widetilde N)_{(x,Z)}))=0,
\quad\forall(x,Z)=Z_x\in T'M$.

iii). $PN(X,QY)=0\ ,\quad\forall X,Y\in{\cal X}(T'M)$.

iv). $PN$ is symmetric with respect to the two lower indices and
it is positive homogeneous of degree 0.
\end{prop}
In the literature on Finsler spaces, one can find a lot of tensor fields having the
properties i)--iv) (see, e.g., \cite{R} and \cite{Ma}) and by following
the known
examples many other tensors of the same kind can be constructed.
From $\widetilde N$ one can obtain many tensors like $N$ by
following methods, which are
analogous to the ones previously used for the construction of the connection
$\nabla^1$. Hence, we omit
them for the sake of brevity.

Now, let us consider a point $x\in M$. Then the ({\it closed}) {\it indicatrix having $0_x$
as its center and $\rho>0$ as its radius} is the set
$\overline{B}_\rho(0_x)=\{(x,X)\in T_xM/F(x,X)\leq\rho\}$ and it is the closure
of the open indicatrix $B_\rho(0_x)$, defined in the Introduction.
$\overline{B}_\rho(0_x)$ and $B_\rho(0_x)$ are both convex, that is for any
$X_x,Y_x$ elements of $\overline{B}_\rho(0_x)$ ($B_\rho(0_x)$) the vector
$tX_x+(1-t)Y_x$ belongs to $\overline{B}_\rho(0_x)$ ($B_\rho(0_x)$), for any
$t\in[0,1]$.

Finally, it is easy to see that:
\begin{prop}
Let $x\in M$, the following assertions are
true:

i). The family $\{B_\rho(0_x)\}_{\rho>0}$ is a basic system of open
neighbourhoods of $0_x$ in $T_xM$.

ii). The sets $\widetilde A=\cup_{z\in A}B_\rho(0_z)$, with $A$ open
neighbourhood of $x$ in $M$ and $\rho>0$ form a basic system of open
neighbourhoods of $0_x$ in $TM$.
\end{prop}

Now, we are in position to prove Proposition 1.1:
\begin{proof}
Let $\nabla$ be a Finslerian Connection and let $x\in M$. Since the Berwald
Connection $\nabla'$ induced by $\nabla$ induces a spray,
by Proposition 2.2, there exists an open neighbourhood $\widetilde W$
of the zero section $\sigma(M)$ in $TM$ and an open subset $W$ of $M$
such that
$(\pi,Exp):\widetilde W\to W\times W$ is a diffeomorphism of class $C^1$, where
$Exp$ is the exponential map of the spray defined by $\nabla'$.
Moreover, ii) of Proposition 2.6 ensures us that there exists a positive real
number $\eta$ and an open neighbourhood $A$ of $x$ in $M$, such that
$\widetilde A=\cup_{z\in A}B_\eta(0_z)\subseteq\widetilde W$. Then, by using ii) of
Proposition 2.3 and i) of Proposition 2.6, there also exists a real number
$\varepsilon$, with $0<\varepsilon\leq\eta$ such that $B_\varepsilon(x)\times
B_\varepsilon(x)\subseteq(\pi,Exp)(\widetilde A)\subseteq W\times W$, where
$B_\varepsilon(x)=Exp_x(B_\varepsilon(0_x))$ is an open neighbourhood of $x$. Finally, i)
follows, since $Exp_x:B_\varepsilon(0_x)\to B_\varepsilon(x)$ is a $C^1$--diffeomorphism.

To prove ii), we consider $y,z\in B_\varepsilon(x)$ with $y\not=z$. It results:
\begin{eqnarray*}
0_y\not=(\pi,Exp)^{-1}(y,z)=X_y^z\in B_\eta(0_y)\subseteq\widetilde W
\ .
\end{eqnarray*}
Then, i) of Proposition 1.1, implies that the curve $c(t)=Exp_y(tX^z_y)$,
with $t\in[0,1]$, is the unique geodesic lying entirely in $W$ and joining
$y$ and $z$. The curve $c$ is also the unique geodesic that joins $y$ and
$z$ lying entirely in $B_\eta(y)$. In fact, we shall see that $c$ has
length lesser than $\eta$.
Let $\beta=(c,\dot c):[0,1]\to T'M$ the complete lift of $c$
to $TM$,
for any $t\in[0,1]$. Then, it results
$\dot \beta=\dot c^ie_i+\ddot c^{\hat i}\varepsilon_{\hat i}$. Hence it follows:
\begin{eqnarray*}
g_{(c(t),\dot c(t))}(\dot c(t),\dot c(t))=
G_{\beta(t)}(\dot\beta(t),\dot\beta(t))\ ,\quad\forall t\in[0,1]
\ .
\end{eqnarray*}
Moreover, being $c$ a path of $\nabla'$, we have:
\begin{eqnarray*}
(\nabla_{\dot\beta}{\dot\beta})(t)=Q_{\beta(t)}(Z_{\beta(t)})\ ,
\quad\forall
t\in[0,1]
\end{eqnarray*}
where $Z$ is a suitable vector field defined on an open neighbourhood of
$\beta([0,1])$ in $T'M$.
From the previous identities and ii) of the previous proposition  we obtain:

\begin{eqnarray*}
& &
\frac{d}{dt}g_{(c(t)\dot c(t)}(\dot c(t),\dot c(t))=
(\nabla_{\dot Y} G)_{\beta(t))}(\dot\beta(t),\dot\beta(t))+
\nonumber
\\
& &
2G_{\beta(t)}(\dot\beta(t),Q_{\beta(t)}(Z_{(\beta(t))}))=0\ ,
\ \quad \forall t\in[0,1]
\ .
\end{eqnarray*}
Hence, we have:
\begin{eqnarray}
F^2(c(t),\dot c(t))=g_{(c(t),\dot c(t))}(\dot c(t),\dot
c(t))=F^2(y,X^z_y)<\eta^2
\ ,\quad\forall t\in[0,1]
\ ;
\end{eqnarray}
that is $c$ has length lesser than $\eta$. Moreover,
for $y=z\in B_\varepsilon(x)$
the unique geodesic joining $y$ and $z$ in $W$ is the constant curve.
Now, we observe that, the mapping $Exp_y:B_\eta(0_y)\to B_\eta(y)$ is a
$C^1$--diffeomorphism, because of ii) of Proposition 2.2, by using the
inclusion $B_\eta(0_y)\subseteq\widetilde W\cap T_yM$.
Finally, we have iii), since
the inclusion
$(\pi,Exp)^{-1}(\{y\}\times B_\eta(x))\subseteq B_\eta(0_y)$ implies
$B_\varepsilon(x)\subseteq B_\eta(y)$.
\end{proof}

Now, we recall that the standard identification $T_XT_xM=T_xM$, induces a
Riemannian Metric on $T'_xM$, for any $X\in T_xM$ and any $x\in M$. By
means of this identification, we prove the corresponding of the Gauss Lemma.
\begin{prop}
Let $x\in M$ and $X\in B_\varepsilon(0_x)$, with $\varepsilon$ as
in Proposition 1.1. Then:

i). It results $F(x,X)=F(Exp_x(X),(DExp_x)_XX)$, where we denoted
by
$(DExp_x)_X:T_XT_xM=T_xM\to T_{Exp_x(X)}M$ the total differential of
$Exp_x$.

ii). For any $Y\in T_xM$, such that $g_{(x,X)}(X,Y)=0$, it
results:
\begin{eqnarray*}
g_{(Exp_x(X),(DExp_x)_XX)}((DExp_x)_XX,(DExp_x)_XY)=0.
\end{eqnarray*}
\end{prop}
\begin{proof}
If $X=0_x$, then the assertion is trivially true. Hence, we suppose
$X\not=0_x$.
\par
In order to prove i) let us denote by $c_X:[0,1]\to M$ the geodesic of
the Berwald connection $\nabla'$, then we have:
\begin{eqnarray*}
(DExp_x)_{tX}X=\dot c_X(t)=\beta(t)\in T'M\ ,\quad\forall t\in[01]\ ;
\end{eqnarray*}
and 
\begin{eqnarray*}
\frac{d}{dt}g_{(c_X(t),\dot c_X(t))}(\dot c_X(t),\dot c_X(t))=0
\end{eqnarray*}
Hence, from Equation (2.2), it follows:
\begin{eqnarray*}
g_{(c_x(t),\dot c_X(t))}(\dot c_X(t),\dot c_X(t))=F^2(x,X)\ ,\quad\forall
t\in[0,1]
\ .
\end{eqnarray*}
Consequently

\begin{eqnarray}
F(c_X(t),\dot c_X(t))&=&F(Exp_x(tX),(DExp_x)_{tX}X)
=
F(x,X),\
\nonumber
\\
& &
\qquad\forall t\in[0,1]\ ;
\end{eqnarray}
and the assertion follows for $t=1$.

Now, we prove ii). Let $Y\in T'_xM$, such that
$g_{(x,X)}(X,Y)=0$. Consider a curve, denoted again by
$X:(-a,a)\to B_\varepsilon(0_x)$, with
$a>0$ and $\varepsilon$ as in Proposition 1.1,
having $X(0)=X$, $\dot X(0)=Y$ and
$F(x,X(s))=F(x,X)=r<\varepsilon$, for each $s\in(-a,a)$. Being $tX(s)\in
B_\varepsilon(0_x)$, we can define the mapping:
\begin{eqnarray*}
\lambda(s,t)=Exp_x(tX(s))=c_{X(s)}(t)\ ,\quad (s,t)\in(-a,a)\times[0,1]
\ .
\end{eqnarray*}
Then:
\begin{eqnarray*}
(\frac{\partial\lambda}{\partial s})_{(s,t)}=(DExp_x)_{tX(s)}t\dot X(s)
\end{eqnarray*}
and
\begin{eqnarray*}
(\frac{\partial\lambda}{\partial t})_{(s,t)}
=(DExp_x)_{tX(s)}X(s)=\dot c_{X(s)}(t)
\ ,
\end{eqnarray*}
for any $(s,t)\in(-a,a)\times[0,1]$.

Hence, it results:
\begin{eqnarray}
F(\lambda(s,t),(\frac{\partial\lambda}{\partial t})_{(s,t)})=
F(x,X(s))=r\ ,\quad\forall
(s,t)\in(-a,a)\times[0,1]\ .
\end{eqnarray}
Now, we set:
\begin{eqnarray*}
\beta(s,t)=(\lambda(s,t),(\frac{\partial\lambda}{\partial t})_{(s,t)})\in T'M,
\quad\forall(s,t)\in(0,a)\times[0,1]
\ .
\end{eqnarray*}
Then, we have:
\begin{eqnarray*}
& &
\bigl(\frac{\partial}{\partial t}
g_{(\lambda,\frac{\partial\lambda}{\partial t})}
(\frac{\partial\lambda}{\partial s},
\frac{\partial\lambda}{\partial t})\bigr)_{\vert(s,t)}=
G_\beta
\bigl(\nabla_{\frac{\partial\beta}{\partial s}}
\frac{\partial\beta}{\partial t},
\frac{\partial\beta}{\partial t})_{\vert(s,t)}=
\nonumber
\\
& &
\frac{1}{2}\bigl(\frac{\partial}{\partial s}G_\beta
(\frac{\partial\beta}{\partial t},
\frac{\partial\beta}{\partial t})\bigr)_{\vert(s,t)}-
\frac{1}{2}\bigl(\nabla_{\frac{\partial\beta}{\partial s}}
G)(\frac{\partial\beta}{\partial t},
\frac{\partial\beta}{\partial t})\bigr)_{\vert(s,t)}
\ .
\end{eqnarray*}
Hence, because of ii) of Proposition 2.3, recalling the definition of $\beta$,
we have:
\begin{eqnarray*}
\frac{\partial}{\partial s}
F^2(\lambda(s,t),(\frac{\partial\lambda}{\partial t})_{(s,t)})=0\ ,\quad
\forall(s,t)\in(-a,a)\times[0,1]\ .
\end{eqnarray*}
Consequently:
\begin{eqnarray}
& &
g_{(Exp_x(tX(s)),(DExp_x)_{tX(s)}X(s))}(
(DExp_x)_{tX(s)}t\dot X(s),
(DExp_x)_{tX(s)}X(s))=0
\nonumber
\\
& &
\forall(s,t)\in(-a,a)\times[0,1]\ .
\end{eqnarray}
Finally, we obtain ii) by putting $t=1$ and $s=0$ in the previous identity.
\end{proof}
\par
For the sequel, we need the following lemma, which is proved in \cite{H}.
\begin{lem}
Let $x\in M$ and $Z\in T'_xM$. Then for each $Y\in T_xM$
such that $g_{(x,Z)}(Z,Y)=0$ it results $F(x,Y+Z)\geq F(x,Z)$ and equality
holds if and only if $Y=0_x$.
\end{lem}
Then, we can prove the Proposition 1.2.
\begin{proof}

First, we prove i).

We can suppose $\widetilde b(s)\not=0_x$, for any $s\in(0,s_1]$. Then, the
function $r(s)=F(x,\widetilde b(s))$ is not zero, for each $s\in(0,s_1]$.
Hence, we can put:
\begin{eqnarray*}
X(s)=\frac{1}{r(s)}\widetilde b(s)\ ,\quad\forall s\in(0,s_1]
\ .
\end{eqnarray*}
With these notations, it is easy to see that:
\begin{eqnarray}
\widetilde b(s)=r(s)X(s)\ ,
\end{eqnarray}

\begin{eqnarray}
g_{(x,X(s))}(X(s),X(s))=1
\end{eqnarray}

\begin{eqnarray}
\dot b(s)=(DExp_x)_{\widetilde b(s)}\dot r(s)X(s)+
(DExp_x)_{\widetilde b(s)}r(s)\dot X(s))\ ,
\end{eqnarray}

\begin{eqnarray}
g_{(x,X(s)}(X(s),\dot X(s))=0\ ,
\end{eqnarray}
for any $s\in(0,s_1]$.

Then, from ii) of Proposition 2.7 and from (2.4), it follows:
\begin{eqnarray}
g_{(b(s),(DExp_x)_{\widetilde b(s)}\widetilde b(s))}
\bigl((DExp_x)_{\widetilde b(s)}\dot
r(s)X(s),(DExp_x)_{\widetilde b(s)}r(s)\dot X(s)\bigr)=0
\end{eqnarray}
for any $s\in(0,s_1]$.
Moreover, i) of Proposition 2.7 and equations (2.4) and (2.5) imply:
\begin{eqnarray}
F^2(b(s),(DExp_x)_{\widetilde b(s)}\dot r(s)X(s))=\dot r^2(s),
\quad\forall s\in(0,s_1]
\ .
\end{eqnarray}
Then, from (2.7), (2.9), (2.10) and from the Lemma 2.1, we obtain:
\begin{eqnarray}
F^2(b(s),\dot b(s))\geq\dot r^2(s)\ ,
\quad\forall s\in(0,s_1]
\ ,
\end{eqnarray}
hence:
\begin{eqnarray*}
L(b)=\int_\rho^{\sigma_1}\vert\dot r(s)\vert ds\geq\vert r(s_1)-r(\rho)\vert\ ,
\quad \forall \rho\in(0,s_1]
\end{eqnarray*}
and i) is true.

Assertion ii) trivially holds.

Under the assumptions of iii), by using the previous notations, we have
$L(c_X)=L(b)=r(s_1)$. Moreover, being i) true, it results
$L(b_s)\geq L(c_{\widetilde b(s)})=r(s)$, with $b_s=b_{\vert[0,s]}$, for any
$s\in(0,s_1]$. Hence, it follows:
\begin{eqnarray*}
\int_\rho^{s_1}(F(b(s),\dot b(s))-\vert\dot r(s)\vert)ds\leq0\ ,
\quad\forall\rho\in(0,s_1]\ .
\end{eqnarray*}
Since (2.12) holds, the previous inequality implies:
\begin{eqnarray*}
\int_\rho^{s_1}(F(b(s),\dot b(s))-\vert\dot r(s)\vert)ds=0\ ,
\quad\forall\rho\in(0,s_1]\
\end{eqnarray*}
and, for continuity reasons, we have
\begin{eqnarray*}
F(b(s),\dot b(s))=\vert\dot r(s)\vert\ ,\quad\forall s\in[0,s_1]\ .
\end{eqnarray*}
Consequently, by using (2.8), (2.10), (2.11) and the previous Lemma we have:
\begin{eqnarray*}
(DExp_x)_{\widetilde b(s)}r(s)\dot X(s)=0\ ,\quad \forall s\in[0,s_1]\ .
\end{eqnarray*}
Then, $\dot X(s)=0$, for any $s\in[0,s_1]$, because $Exp_x$ has maximal rank
along the curve $\widetilde b$ and $r(s)\not=0$ for any $s\in(0,s_1]$.
Hence, iii) is true.
\end{proof}

From the above proposition, it trivially follows:
\begin{prop}
Let $x\in M$. Moreover, let us consider
$\varepsilon=\varepsilon(x)$ and
$\eta=\eta(x)$ as in Proposition 1.1. Then, a geodesic of
length lesser than $\eta$ starting
from an arbitrary point $y\in B_\varepsilon(x)$
is a curve of minimal length between its end points.
\end{prop}
Let $x\in M$ and $B_\rho(0_x)\subseteq\widetilde W\cap T_xM$. Moreover,
let $u=(u_i)_{1\leq i\leq n}$ be a basis of $T_xM$. Then, $u$ defines a
mapping, which by an abuse of notations, we denote again by $u:{\mathbb R}^n\to
T_xM$ obtained by putting $u(\xi)=\xi^iu_i$, for any $\xi=(\xi^i)_{1\leq
i\leq n}\in{\mathbb R}^n$. Furthermore, the mapping
$\psi=u^{-1}\circ Exp_x^{-1}:B_\rho(x)\to u^{-1}(B_\rho(0_x))$ is a
$C^\infty$--diffeomorphism on $B_\rho(x)-\{x\}$ and $C^1$--differentiable on
$x$. Consequently, $(B_\rho(x),\psi)$ is a chart of the
$C^1$--differentiable manifold structure canonically induced on $M$ by the
considered
structure of $C^\infty$--differentiable manifold on $M$.

Then, we can prove:
\begin{prop}
Let $x\in M$ and $\varepsilon=\varepsilon(x)$ as in Proposition
1.1. There exists $\varepsilon_0=\varepsilon_0(x)\in(0,\eta)$
such that for any non constant
geodesic $c:[0,a]\to M$, with $a>0$, satisfying the conditions
$c(t_0)\in\partial B_\varepsilon((x)=Exp_x(\partial B_\varepsilon(0_x))$ and
$\dot
c(t_0)\in T\partial B_\varepsilon(x)$, for some $t_0\in(0,a)$ and
$\varepsilon\in(0,\varepsilon_0)$
and for each $\mu\in(0,1)$ there exists $\rho=\rho(\mu,\varepsilon)$ with the
following property:
\begin{eqnarray*}
d(x,c(t))\geq d(x,c(t_0))+\mu(t-t_0)^2\ ,\quad \forall
t\in[t_0-\rho,t_0+\rho]\ ,
\end{eqnarray*}
where $d$ is the Finslerian distance function.
\end{prop}
\begin{proof}
We fix a basis of $T_xM$ and consider the $C^1$--differentiable chart\\
$(B_\eta(x),\psi)$, with $\psi=u^{-1}\circ Exp_x^{-1}$. Moreover, we can
suppose $x\not=c(t)\in B_\eta(x)$, for any $t\in(0,a)$. Then all the
derivatives of the function $f(t)=d^2(x,c(t))=F^2(x,\gamma(t))$, being
$\gamma(t)=\psi(c(t))$, for any $t\in(0,a)$, are defined. Hence by
using the homogeneity conditions, we get:
\begin{eqnarray*}
& &
f(t_0)=\varepsilon^2\ ;\quad \dot f(t_0)=
g_{(x,\gamma(t_0))}(\gamma(t_0),\dot\gamma(t_0))\
;
\nonumber
\\
& &
\ddot f(t_0)=
g_{(x,\gamma(t_0))}(\dot\gamma(t_0),\dot\gamma(t_0))+
g_{(x,\gamma(t_0))}(\gamma(t_0),\ddot\gamma(t_0))\ .
\end{eqnarray*}
Therefore, the proof follows as in \cite{Kl} (cf. \cite{H}, too).
\end{proof}

The proof of Proposition 1.3 follows from the previous proposition in a
trivial way.

We can also omit the proof of the following Lemma, because it needs only the
compactness of boundary of convex neighbourhoods, which holds, because an
indicatrix is always compact and the exponential map is a diffeomorphism
near to any point. Hence, its proof follows as in \cite{Kl}.
\begin{lem}
Let $p\in M$ and $\rho>0$ be a positive real number such
that  $Exp_p$ is defined on the ball $B_\rho(0_p)\subseteq T_pM$. Then,
every $q\in M$, with $d(p,q)<\rho$, can be joined to $p$ by a minimizing
geodesic.
\end{lem}

Finally, the proof of Proposition 1.4 is the same as in the Riemannian case
(see, e.g., \cite{Kl}).

\end{document}